# A surrogate-based approach to accelerate the design and build phases of reinforced concrete bridges

Mouhammed ACHHAB, Pierre JEHEL, Fabrice GATUINGT

*Université Paris-Saclay, CentraleSupélec, ENS Paris-Saclay, CNRS, LMPS – Laboratoire de Mécanique Paris-Saclay, 91190, Gif-sur-Yvette, France,*

**ABSTRACT** Integrating uncertainties in the design process of reinforced concrete rail bridges, in a fully probabilistic framework, makes their design more complex and challenging. To propagate these uncertainties and convey their influence on the performance of the engineering system, a high-dimensional design space is supposed to be explored. A great challenge to be considered here lies in the computational burden as conducting such an exploration campaign requires substantial calls to computationally expensive finite element simulations. To address this challenge, a surrogate model mapping the design space to the reinforced concrete bridge performance functions is developed in the context of an active learning algorithm. The importance of this model lies in its ability to explore as many design scenarios as possible with minimal computational resources and classify the design scenarios into failure and safe scenarios. This work considers a 4-span reinforced concrete bridge deck. A multi-fiber finite element model of this beam is developed in Cast3m to generate the required design of experiments for the surrogate model. A performance comparison is undertaken to evaluate the Kriging surrogate model effectiveness with and without active learning while the reliability of Kriging predictions is also assessed in comparison to PC-Kriging.

## I. INTRODUCTION

The growing demand for rail bridges is driven by the need to enhance connectivity, promote fast and effective transportation systems, and address the challenges of urban densification. Among all the rail bridge types, reinforced concrete bridges are a common choice in modern structural design. The scale of rail bridges, their complex operating conditions, and the concrete non-linear mechanical behavior make the design of such structures complex and time demanding.

For the design and analysis process, bridge engineers have long relied on numerical models using finite element analysis tools along with design standards and specifications as the Eurocodes in Europe. The limit states principle is a key concept in the Eurocode design rules. Limit states split into two main groups: ultimate limit states to guarantee the safety of people and structure, and serviceability limit states to ensure that the comfort of the users is satisfied, and the appearance of the bridge is not excessively affected. In the context of bridge engineering, ultimate limit states can be represented by flexural or shear capacities, while serviceability limit states include limits on deflection, cracking widths, and vibrations.

In the Eurocodes, loads acting upon the structure, as well as structural capacity, are corrected by security factors. Such factors were introduced in design codes to account for some types of uncertainties, including loads, material properties, and construction and manufacturing tolerances.

However, unexpected situations can happen in the design and build process that lack reliable statistical or empirical representation, making them difficult to quantify using predefined factors. For instance, in the design and build process, previously undetected underground structures or utility networks may obstruct foundation placement, preventing the execution of the design as originally planned.

This paper is about proposing an efficient reinforced concrete bridge design approach enhancing structural adaptability and minimizing redesign costs and construction delays through the anticipation of unexpected events from the early design stage. These latter can disrupt the continuity between the design and execution phases. Also, engineers are frequently required to modify design parameters, leading to iterative redesigns to accommodate unforeseen conditions. The work presented in this paper focuses on reinforced concrete bridges. It considers that unexpected situations result in geometrical uncertainties, thus making it possible to address our problem using probabilistic tools. This goes beyond the scope of design codes.

We aim at exploring a potentially large space of design scenarios. The main challenge in our exploration is the computational overhead, stemming from the need for numerous, potentially time-consuming finite element simulations. To overcome this, we are developing a surrogate model, as a predictive model, that leverages a limited number of FE simulations for efficient uncertainty propagation. We developed a simulator of a reinforced concrete continuous bridge deck to provide the required simulations for building a surrogate model. This predictive model maps the design space, where uncertainties are modeled as random variables, onto a limit state function of the bridge. In this work, the simulator is an FE model with multifiber beams; the uncertainties considered are the positions of the piers; the limit state is the maximum deflection of the deck.

Our objective is to explore as many design scenarios as possible with minimal computational resources and classify these design scenarios into failure or safe scenarios. We use the Kriging approach to develop our surrogate model, as one of the most widely and commonly used approaches in the field of computational experiments, proving its efficiency and performance in many domains, including structural engineering. We evaluate in this paper its performance with the variation of the number of simulations used to generate the metamodel. As we aim to classify our design scenarios, we should be able to accurately approximate the limit state function (LSF). Besides, the number of costly finite element (FE) simulations should be reduced, which means the surrogate model must reliably detect and achieve high precision near the failure boundary—the critical region separating safe and failure domains. By strategically selecting new sample points in these regions, active learning enhances model efficiency, ensuring reliable predictions with minimal computational cost. To assess the performance of the active learning algorithm, we conduct a study to check the influence of combining this approach with Kriging in our problem. Also, within the active learning framework, we compare the performance of traditional universal Kriging with polynomial function regression against the performance of the PC-Kriging surrogate model.

In Section II, the reinforced concrete bridge model used to develop the surrogate is presented. Section III introduces the surrogate modeling strategy implemented for this study. In section IV, numerical applications are shown, and the obtained results are discussed. Section V is devoted to the conclusions and ongoing research work.

## II.  BRIDGE DECK FEM

A. *Model description*

We consider a reinforced concrete bridge deck. A finite element model of a four-span continuous reinforced concrete deck is developed using Cast3M software. The beam is 40 m long, supported on five supports, preventing displacements in all directions. Four concentrated loads are applied to the beam in different locations. The model of the beam representing the applied loads, the supports, and the geometrical parameters is shown in FIGURE 1.

B. *Parametric FE model*

A multi-fiber (Combescure, 2007) finite element model was considered in our problem. Known for its cost-effectiveness, this approach has been implemented in Castem 2000 software in 1994. This method is of interest when no fine mesh is needed to detect important local information. It was developed as a compromise between the simplicity of the beam element and the accuracy of the expensive 3D nonlinear models. It consists in applying the beam theory (Euler- Bernoulli (Öchsner, 2021a) or kinematic Timoshenko (Öchsner, 2021b) along the longitudinal axis of the structure. Each Gauss point of the beam longitudinal elements represents a transverse section. At the level of the sections, multiple fibers are considered for a finer discretization. In a reinforced concrete FE model, fibers are categorized into three groups, each representing different material behavior. One group models unconfined concrete, corresponding to the concrete cover. Another models confined concrete, which accounts for the increased strength due to stirrups. The longitudinal steel is represented by attributing its mechanical behavior to a third group.

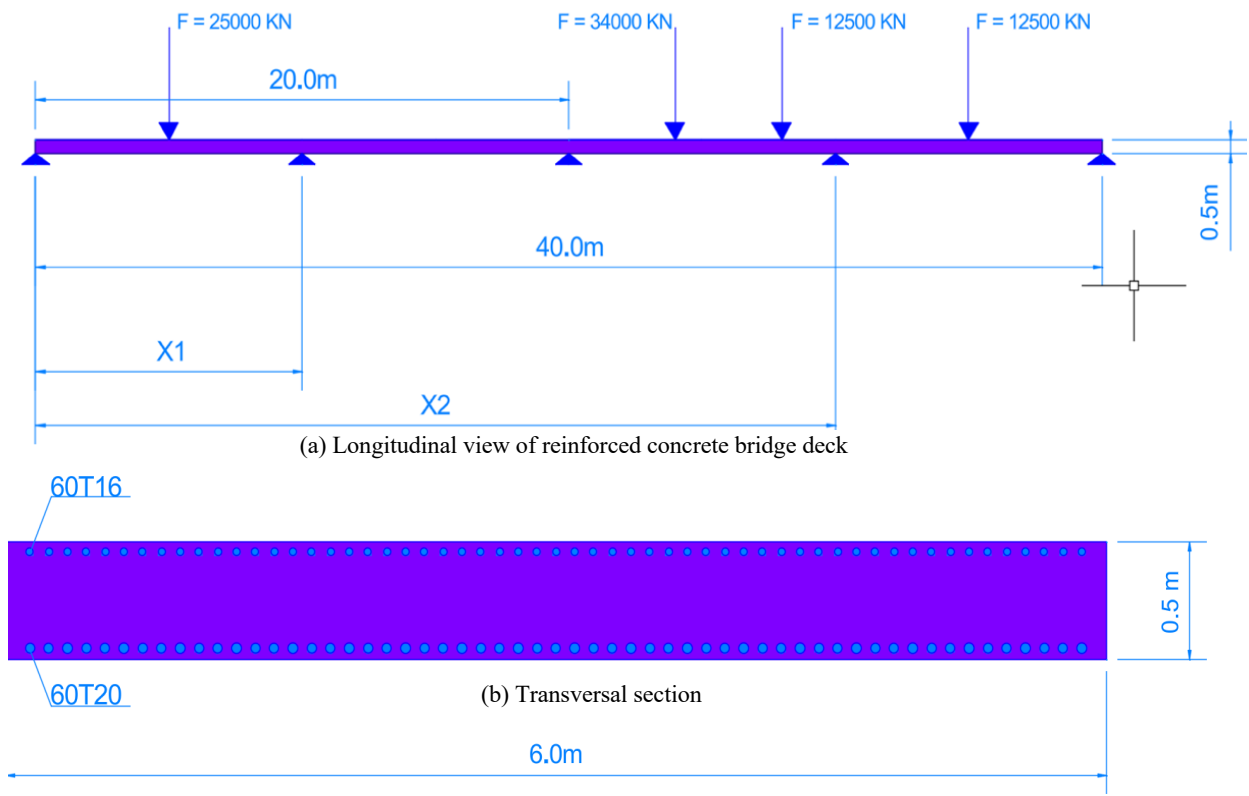

(a) Longitudinal view of reinforced concrete bridge deck

(b) Transversal section

FIGURE 1. **Model of the problem**

As shown in FIGURE 2, the beam was discretized longitudinally using the Timoshenko element. At the level of the section, steel fibers were modeled as point elements POJS with one integration point by rebar, and the concrete was modeled as a QUAS element with four integration points by fiber. Eurocode model (Beton BAEL in Cast3M) (European Committee for Standardization (CEN), 2004) and kinematic plastic hardening (Parfait Uni in Cast3M) (McDowell, 1987) were used for concrete and steel, respectively. A parametric FE model (FEM) was developed, where two parameters are treated as random variables to account for the variability in the positions of the second and fourth piers. These parameters are represented by the random vector $\mathbf{X} = [X_1, X_2]^T$, where $X_1$ and $X_2$ are uniformly distributed between 3 and 18 m, and 23 and 38 m, respectively. The position of the first, middle, and fifth piers are fixed at 0m, 20m, and 40m, respectively. With this setting, we intend to investigate the impact of modifying piers' positions on the maximum deflection in the bridge deck.

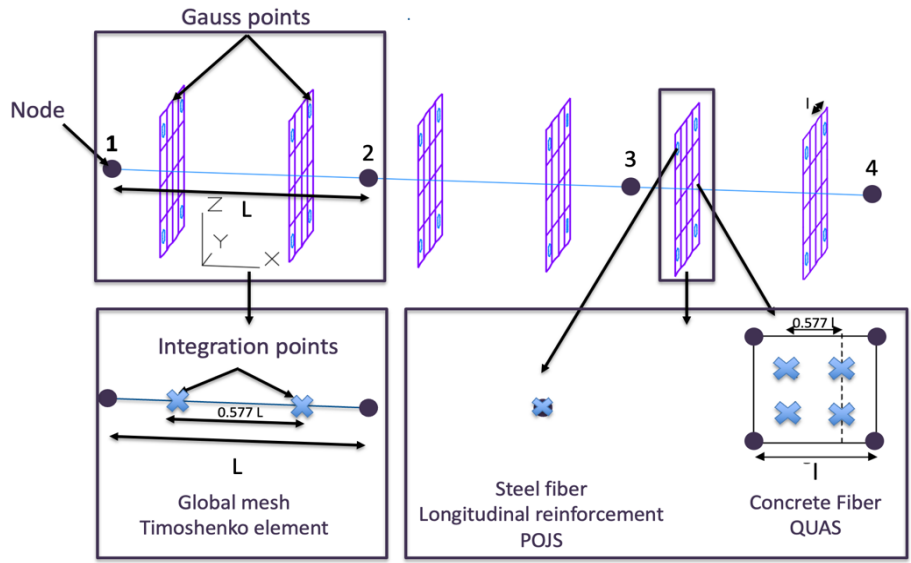

**FIGURE 2.** Multifiber element

C. High-fidelity limit-state model

The deflection limit state function is represented by the function $g(\mathbf{x}) = L - q(\mathbf{x})$. L represents the limit the deflection should not reach for the bridge to stay in service conditions, and $q(\mathbf{x})$ is the maximum deflection calculated using the FE model. To approximate the limit state function with accuracy, we have done 10,000 evaluations of $g$ for different values of $\mathbf{x}$ selected using the Latin Hypercube sampling technique (Stein, 1987). A computational pipeline was developed using MATLAB and Cast3M software. Latin Hypercube Sampling (LHS) is implemented in MATLAB to select the desired values of the input variables $\mathbf{x}$. For each selected $\mathbf{x}$, the Cast3M software is called automatically by MATLAB to compute the corresponding maximum deflection value. All the results of Cast3M calculations, along with their associated $\mathbf{x}$ values, are stored back in MATLAB, where the post-processing is performed. The plot of the 10,000 evaluations is shown in FIGURE 3, where the limit between the failure (red) zone and the safe zone (blue) is well defined. Based on these 10,000 evaluations, the reference structural failure probability (Der Kiureghian, 1996) $P_{f(ref)} = P(g(\mathbf{x}) \leq 0) = \frac{N_f}{10,000}$ is calculated ($N_f$ represents the number of failure points with $g(\mathbf{x}) \leq 0$). The obtained reference probability of failure is $P_{f(ref)} = 0.0671$.

### III.     SURROGATE MODEL OF THE LIMIT STATE FUNCTION

In the previous section, we used 10,000 simulations for the modeling process, while in real life, one simulation of a bridge is computationally expensive. Efficient computational strategies should be used to approximate the function $g(\boldsymbol{x}) = 0$ that splits the design scenarios into admissible or not. Our approach is based on the construction of a surrogate or metamodel of $g(\boldsymbol{x})$. To generate the necessary data for building a surrogate model, n samples $\boldsymbol{S} = \{\boldsymbol{x}^{(1)}, ..., \boldsymbol{x}^{(n)}\}$ were selected. LHS was used with uniformly distributed X1 ~ U (3,18) and X2 ~ U (23,38) to ensure a well-distributed set of input samples across the specified ranges for both pier positions. The output $g^{(i)}(\boldsymbol{x}) = g(\boldsymbol{x}^{(i)})$ for each sample point is calculated. The resulting input-output data set $\boldsymbol{S}$ and $\mathcal{G} = [g(\boldsymbol{x}^{(1)}) ... g(\boldsymbol{x}^{(n)})]^T$ form the foundation for the surrogate model development. Our goals in this context are:

- To develop a surrogate model using kriging with different numbers of samples and evaluate their performance.

- To develop an active learning algorithm with Kriging to evaluate its influence

- To develop an active learning algorithm with Polynomial Chaos Kriging (PC-Kriging: a combination of both Polynomial Chaos expansion and (PCE) and Kriging)

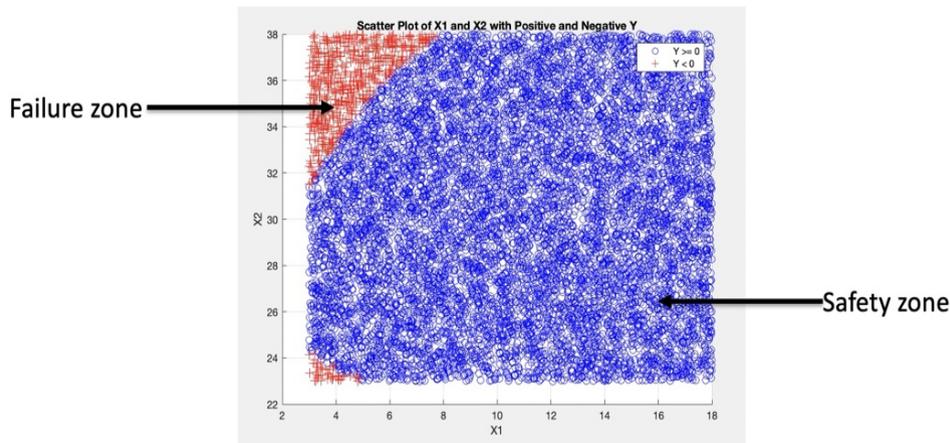

**FIGURE 3.** Results of the reference problem using FEM simulations

#### A.   Kriging

Gaussian processes (Rasmussen and Williams, 2008), known as Kriging surrogate models, were first introduced for geostatistical applications (Krige, 1951). Since then, they have evolved into one of the most effective methods for approximating computer experiments. Kriging surrogate models have served as performant emulators for computationally expensive finite element analysis problems in different contexts such as design optimization (Pellegrino et al., 2015), reliability analysis (Dubourg et al., 2013), and computer model calibration (Kennedy and O'Hagan, 2001). They offer a flexible and robust approach to approximating complex functions, providing accurate predictions and reliable uncertainty estimates even in scenarios with limited or noisy data.

Having a set of (observed) data points $S$ and their corresponding solutions (observations) $\mathcal{G}$, the goal of the kriging algorithm is to allow for efficient prediction of $g(\boldsymbol{x})$ at untried inputs. With Kriging, the function $g(\boldsymbol{x})$ is modeled by $\hat{g}(\boldsymbol{x})$, as a realization of a Gaussian process.  The model is defined as a sum of a deterministic part (the prior mean function $m(x)$ and a Gaussian process with zero mean and covariance function $C(x, x')$ as stated in Eq.(1):

$$\hat{g}(x) = m(x) + Z(x) = \beta^T f(x) + Z(x) \qquad (1)$$

where $Z(x) \sim \mathcal{GP}(0, C(x, x'))$ represents the random spatial fluctuations of $\hat{g}(x)$ around the mean. The prior mean or the trend $m(x)$ reflects our belief about the underlying function, while the covariance $C(x, x') = \sigma^2 R(x, x', \theta)$ models the spatial correlation between different locations ($\sigma^2$ is the variance of the Gaussian process and $R(x, x', \theta)$ the correlation function between two different points denoted $x$ and $x'$, with correlation length $\theta$.). Both initial mean and covariance are termed "prior" because they encapsulate assumptions about $g(x)$ prior to any observations. Once simulations are made, the posterior distribution of $\hat{g}(x)$ is updated, providing a refined estimate of $g(x)$ as well as the associated uncertainty. The posterior distribution of $\hat{g}(x)$ is stated in (2):

$$\hat{g}(x)|\mathcal{G} \sim \mathcal{N}\left(\mu_{\hat{g}}(x), \sigma^2_{\hat{g}}(x)\right) \qquad (2)$$

where $\mu_{\hat{g}}(x)$ is the posterior mean or Kriging predictor (approximation of $g(x)$ at unknown points ($x$) and $\sigma^2_{\hat{g}}(x)$ its posterior variance. More details about the Kriging theory and the estimation of $\mu_{\hat{g}}(x)$ and $\sigma^2_{\hat{g}}(x)$ can be found in (Sacks et al., 1989);(Lataniotis et al., 2018)

The regression functions $f(x)$ used in our problem are 4th-order polynomial functions. The choice is thought to well-fit the approximated function as the deflection in a simply supported beam depends on the 4th order of the span length). The hyperparameters $\sigma^2$, $\theta$, and $\beta$ coefficients are estimated using maximum likelihood of the observations $\mathcal{G}$. The considered correlation function is the well-known Matern function (Petit et al., 2022).

B. *Polynomial Chaos Kriging (PC-Kriging)*

This method a combination of two popular existing methods, kriging and Polynomial Chaos Expansion (first mentioned and introduced in (Wiener, 1938)). The idea of this method and the proof of its performance were shown in a paper in 2015 (Schoebi et al., 2015). This method works as a special case of the universal kriging algorithm, where the trend is replaced by a set of orthonormal polynomials. The performance of this method lies, as stated by Eq.(3), in the combination of both a global PCE approximation and a local Kriging approximation:

$$\hat{g}(x) = \sum_{\alpha \in \mathcal{A}} \xi_\alpha \Psi_\alpha(X) + \sigma^2 Z(x) \qquad (3)$$

$\sum_{\alpha \in \mathcal{A}} \xi_\alpha \Psi_\alpha(X)$ represents the trend of the PC-Kriging developed using polynomial Chaos Expansion (PCE). The type of the polynomials in the trend depends on the distribution functions of the components of the random vector (Xiu and Karniadakis, 2002).

Considering the uniform distributions of each of the random variables $X_1$ and $X_2$, the Legendre polynomials are used.

C. *Active learning surrogate model*

Trying to approximate the limit state function with a minimum number of simulations, in this section, we investigate the potential of an active learning surrogate model. It consists of starting with an initial sampling to perform later an adaptive sampling promoting the selection in the areas where the prediction error of the surrogate is high. The necessary steps used, as shown in FIGURE 4, are the following:

- The initial sampling, referred to as the initial design of experiments, is conducted using Latin hypercube Sampling. For our two-dimensional problem, we determined the number of initial samples using the empirical formula n = max(2*m, 10), where n represents the number of samples and m is the problem's dimension (Moustapha et al., 2022). So, with the Kriging surrogate, 10 points were used in our problem.

- The surrogate model is then developed, as shown in the previous sections

- The active learning algorithm is applied: A learning function is a method to enrich the surrogate model by selecting data points where the model is less accurate to improve its efficiency and accuracy. Different learning functions are presented in the literature Expected Feasibility Function (EFF) (Lv et al., 2015), Expected Improvement (EI) (Emmerich et al., 2008), Fraction of Bootstrap replicate (FBR) (Marelli and Sudret, 2016), and U (deviation number) (Echard et al., 2011). In this work, the deviation number is used among others based on (Moustapha et al., 2022). The deviation number is a function of the kriging variance, and it is based on the theory of misclassification probability. Its formula is shown hereafter: $U(x) = \frac{|\mu_{\hat{g}}|}{\sigma_{\hat{g}}}$ where $\mu_{\hat{g}}$ the mean of kriging and $\sigma_{\hat{g}}$ its standard deviation. The enrichment points are selected so that they minimize the value of the deviation number: $x_{\text{next}} = \text{argmin}(U(x))$

- Finally a convergence criterion is used to stop the iterations. We used the stability of failure probability as stated in Eq.(4)

$$\frac{\left|p_f^{(i)} - p_f^{(i-1)}\right|}{P_f^{(i)}} \leq \epsilon_{P_f} \tag{4}$$

$p_f^{(i)}$ represents the failure probability estimated by Subset simulation (Au and Beck, 2001) at $i^{\text{th}}$ iteration, and $\epsilon_{P_f}$ is the threshold or tolerance. The algorithm converges if the condition in Eq. 4 is verified for three concecutive iterations.

## IV. RESULTS

In this section, we present the numerical implementation and the results derived from the application of our methodology. The objective is to approximate the limit state function of the bridge, where the limit state is defined by the maximum allowable deflection. Specifically, the considered deflection limit is set at $L = x_1/400$. To approximate the limit state function and classify the pier positions as acceptable or not with respect to this deflection limit, we employed and compared the performance of three distinct approaches: Kriging, Active Learning Kriging, and Active Learning PC-Kriging. The high-fidelity model developed in section II.C using 10000 simulations was used to assess the performance of each surrogate model qualitatively and quantitatively. Qualitatively, by visually evaluating the accuracy in reproducing the limit state $g(x) = 0$ *and* representing both the upper and lower failure zones as shown in FIGURE 3,. As well as, by calculating, quantitatively, the relative error as stated in Eq.(5)

$$E_{\text{Pf}} = \frac{|P_{f(ref)} - P_{\text{f}}|}{P_{f(ref)}} \tag{5}$$

where $P_{f(ref)} = 0.0671$ is the reference failure probability obtained from the high fidelity model in section II.C, and $P_{\text{f}}$ is the failure probability estimated by the surrogate model. To compare the performance of the implemented approaches, the evolution of this error metric with the number of numerical simulations is plotted in FIGURE 8.

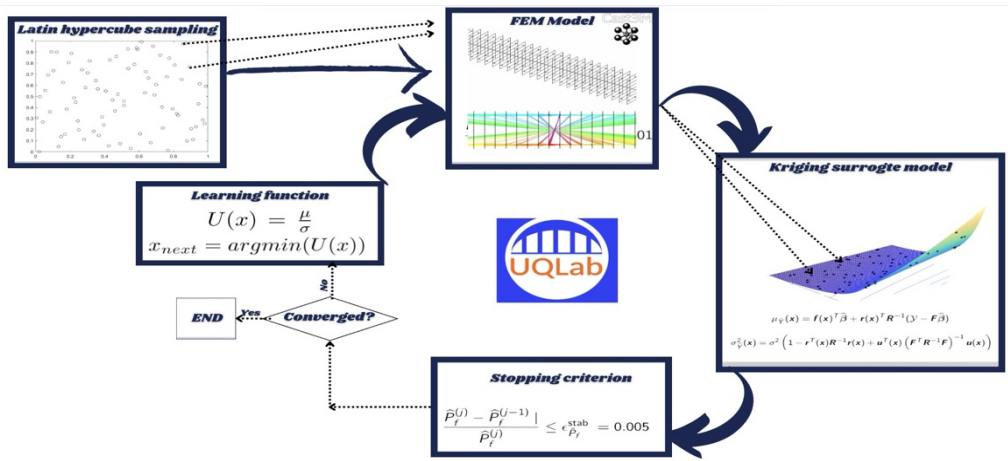

**FIGURE 4.** Active learning steps

For the **Kriging** approach, we developed many surrogate models with different number of simulations (data points), and we evaluated the error $E_{\text{Pf}}$ for each of these surrogates. We started with 10 initial simulations and gradually increased the sample size to 25, 40, 50, 60, 70, and 90 simulations. The results show that $E_{\text{Pf}}$ decreases as shown in FIGURE 8, and the reconstruction of the reference model of FIGURE 3 improves as shown in FIGURE 5, as the number of simulations increases. The surrogate model reached its maximum accuracy with 90 simulations, with $E_{\text{Pf}} = 0.7\%$.

For **Active Learning Kriging**, the convergence of the model was controlled by a convergence threshold $\epsilon_{P_f}$ (see Eq. (4)). Additionally, a maximum simulation cap of 90 was considered, in compliance with the results from the Kriging model. We started with 10 initial finite element simulations and, at each iteration, added one extra point near the limit state $g(x) = 0$ to improve the surrogate model's performance in this critical region. Three experiments were conducted to evaluate Active Learning Kriging with different convergence criteria:

In **experiment 1**, the model was allowed to converge with a relatively loose convergence threshold $\epsilon_{P_f} = 0.005$. The model converged after 24 simulations, with an acceptable error of $E_{\text{Pf}} = 0.68\%$ (See FIGURE 8). FIGURE 6.a shows visually how this method reproduces the results of FIGURE 3.

In **experiment 2**, a more stringent tolerance of $\epsilon_{P_f} = 0.0005$ was used. In this case, the model required 56 simulations to converge, achieving very high accuracy. Again, the error value $E_{\text{Pf}}$ at the convergence point is highlighted in the plot of FIGURE 8. The visual accuracy is high as shown in FIGURE 6.b.

In **experiment 3**: To compare Active Learning Kriging against Kriging, we let the algorithm run until 90 simulations. The plot for this experiment shows the progression of $E_{Pf}$ until 90 simulations.

For **Active Learning PC-Kriging**, we followed a similar approach. We started with 10 simulations and conducted three experiments:

In **experiment 1**: With a convergence criterion of $\epsilon_{P_f} = 0.005$ the model converged after 15 simulations, but it did not accurately the lower failure zone as shown in FIGURE 7.a. The quantitative error $E_{Pf}$ is shown in the plot of FIGURE 8.

In **experiment 2** With $\epsilon_{P_f} = 0.0005$, the model gave the same results as the first experiment.

In **experiment 3**: we let the algorithm run again until 90 simulations. FIGURE 8 shows that the model required 65 simulations to start giving an accurate estimation of the failure probability. Visual evaluation is illustrated in FIGURE 7.b

These results demonstrate that Active Learning Kriging outperforms both Kriging and PC-Kriging active learning in terms of accuracy with fewer simulations.

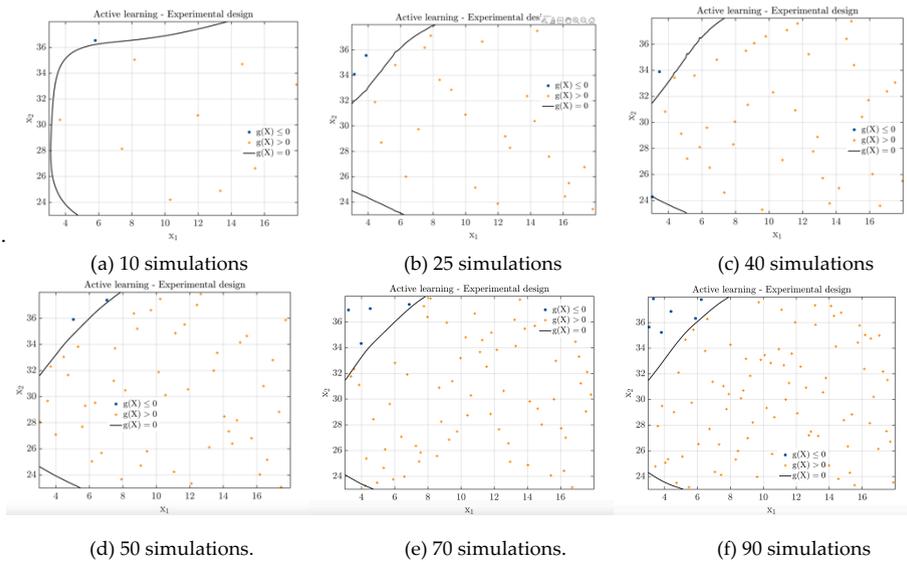

(a) 10 simulations  (b) 25 simulations  (c) 40 simulations

(d) 50 simulations.  (e) 70 simulations.  (f) 90 simulations

**FIGURE 5.** Kriging results (Compare with FIGURE 3.)

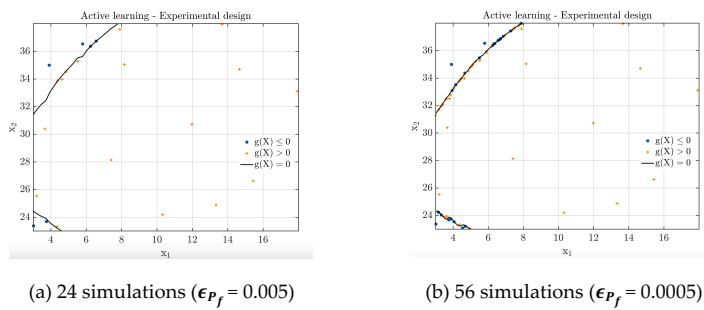

(a) 24 simulations ($\epsilon_{P_f}$ = 0.005)  (b) 56 simulations ($\epsilon_{P_f}$ = 0.0005)

**FIGURE 6.** Active learning Kriging results (Compare with FIGURE 3.)

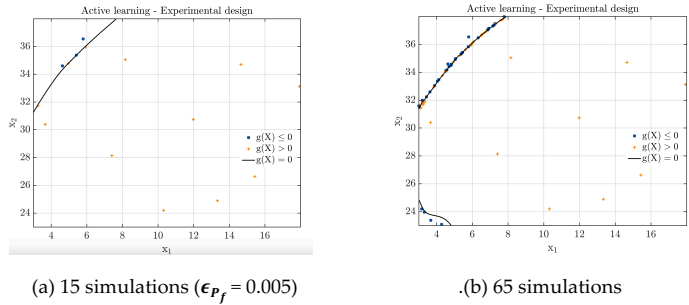

(a) 15 simulations ($\epsilon_{P_f}$ = 0.005)    .(b) 65 simulations

**FIGURE 7.** PC-Kriging active learning results (Compare with FIGURE 3.)

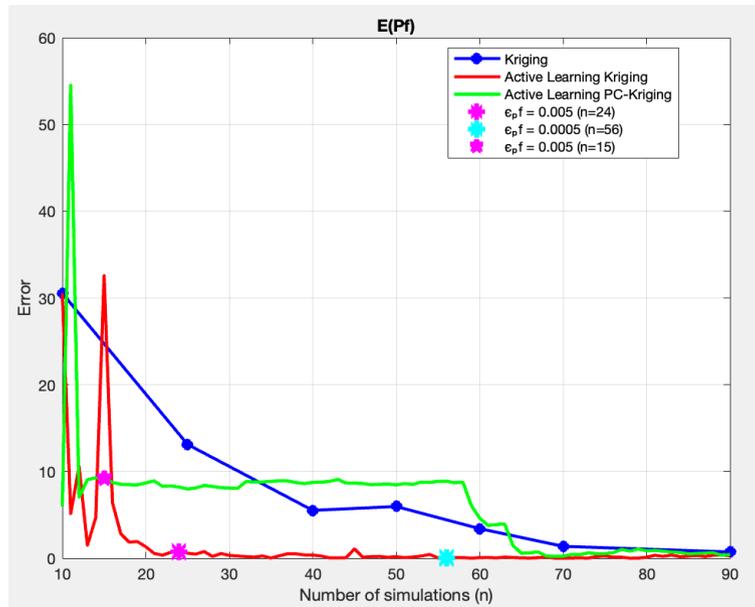

**FIGURE 8.** Comparison of the probability estimation error

## V.    CONCLUSIONS

This study demonstrates the significant potential of surrogate modeling, particularly in accelerating the design and construction processes of reinforced concrete bridges, and increasing the continuity between these two phases. Without the use of surrogate models, approximating the limit state and failure probability would have required over 1000 high-fidelity simulations (e.g., FEM runs), which would be computationally prohibitive. However, by employing surrogate models, we were able to obtain accurate approximations of the limit state function and failure probability using significantly fewer simulations. This not only saved computational resources but also enabled a more efficient exploration of the design space. Among the three methodologies considered—traditional Kriging, Active Learning Kriging, and Active learning PC-Kriging—Active Learning Kriging proved to be the most efficient and accurate. Active Learning Kriging focused simulations on regions with higher uncertainty near the limit state, reducing the number of required simulations. While Kriging needed 90 simulations to achieve a modest reduction in error, Active Learning Kriging converged with just 24 simulations, achieving an error of 0.68%. The PC-Kriging model required more simulations (65) to accurately capture the lower failure zone. While this approach has been very effective in this problem, its full potential becomes more evident in cases involving more complex situations using expensive finite element models. In such cases, surrogate modeling can be invaluable for significantly reducing computational effort, making it a crucial tool for large-scale and computationally expensive problems. The results emphasize the practical advantages of Active Learning Kriging for surrogate modeling, demonstrating that it can reduce the computational burden without sacrificing accuracy, particularly in structural reliability analysis. For the ongoing work, we are considering a mechanical model with different non-linear mechanical models. Also, we are considering a multidimensional problem, including more design parameters in the modelling process. A time-dependent surrogate model could also be developed to consider the influence of materials aging in our problem.

## ACKNOWLEDGEMENT

"The MINERVE project has been financed by the French government within the framework of France 2030."